\documentclass[12pt]{article}

\def\note#1{}

\begin{document}
\baselineskip 21pt
~\vspace{1cm}
\begin{center}

{\Large\bf Non-commutative 4-spheres based on all Podle\'s 2-spheres and
beyond}

\vskip 1cm

{\large {\bf Tomasz Brzezi\'nski$^{1,2}$ and Cezary Gonera$^3$}}

\vskip 0.5 cm
${}^1$Department of Mathematics, University of Wales Swansea\\
Singleton Park, Swansea SA2 8PP, U.K.\\ E-mail:
T.Brzezinski@swansea.ac.uk\\ ~\\
$^{2}$Department of Theoretical Physics, University of 
\L\'od\'z,\\ ul. Pomorska 149/153, 90-236 {\L}\'od\'z,
Poland.\\~\\
${}^3$Department of Theoretical Physics II, University of 
\L\'od\'z,\\ ul. Pomorska 149/153, 90-236 {\L}\'od\'z,
Poland.

\end{center}
\vspace{1 cm}
\begin{abstract}
\begin{quote}
\noindent A wide class of noncommutative spaces, including 4-spheres 
based on all the quantum 2-spheres and suspensions of matrix quantum 
groups is described. For each such space a noncommutative vector 
bundle is constructed. This generalises and clarifies 
various recent constructions of 
noncommutative 4-spheres.
\end{quote} 
\end{abstract}
\thispagestyle{empty}
\newpage 

\noindent {\bf 1.} Recently there has been an upsurge of activity in
constructing noncommutative 4-spheres and corresponding 
projective modules or (instanton) noncommutative vector bundles. 
This has been initiated by a paper by Connes and
Landi \cite{ConLan:man} where, among others, a family of isospectral 
noncommutative
4-spheres labelled by a deformation parameter of unit modulus has been
constructed. Another family of $q$-deformed noncommutative 
4-spheres labelled by a real
deformation parameter $q$ was defined via a suspension of the 
quantum group 
$SU_q(2)$ in
\cite{DabLan:ins}. These
have been followed up in \cite{Sit:mor}, where a
two-parameter family of noncommutative 4-spheres that include the
Connes-Landi one as a special case, was introduced as a noncommutative
two dimensional suspension of one of Podle\'s 2-spheres
\cite{Pod:sph}. In a different direction, motivated by 
the classical Hopf fibering $S^{7}\to S^{4}$, a family of 
noncommutative 4-spheres was 
defined in \cite{BonCic:ins} as coinvariants of the coaction of 
$SU_{q}(2)$ (viewed as a coisotropic subgroup of $U_{q}(4)$) 
on the quantum 7-sphere, i.e.,  as a quantum 
quotient space in the 
sense of \cite{Brz:hom}. 
In all the above cases, projective modules were
constructed with a projector given in terms of a $4\times 4$-matrix with
entries from the noncommutative space, and the corresponding Chern-Connes
character was computed.

The aim of the present note is twofold. First we would like to
generalise the construction in \cite{Sit:mor} by describing a class of 
noncommutative 4-spheres obtained by a two-dimensional noncommutative
 suspension of all the Podle\'s
2-spheres. Second we would like to propose a general construction of
noncommutative manifolds based on the standard formulation of the matrix
quantum group $GL_q(n)$ (the FRT-construction). The q-deformed 4-sphere
in \cite{DabLan:ins}
can then be viewed as a special case of this construction with $n=2$.
\medskip 

\noindent {\bf 2.}
Algebra of functions on the quantum sphere $A(S_{q,s}^2)$ \cite{Pod:sph}
 is defined as a
polynomial algebra generated by $1,x,y,z$ subject to 
the relations
\begin{equation}
zx = q^2xz,\quad yz = q^2zy, \quad xy =
(z -1)(z +s^2), \quad yx = (q^{2}z -1)(q^{2}z
+s^2), 
\label{rel.s2}
\end{equation}
where $q, s$ are complex parameters, $q\neq 0$, $s^{2}\neq -1$. A commutative algebra of functions
on the 2-sphere can be identified with $A(S_{q,s}^2)$ with $q=1$ and $s=0$.
From a purely algebraic point of
view there is no restriction on the values of deformation parameters $s$
and $q$. In all cases $A(S_{q,s}^{2})$ is a subalgebra of the algebra of
functions on the quantum group $SL_q(2)$. If $q^{2}$ and $s^2$ are real then
$A(S_{q,s}^{2})$ is a $*$-algebra with involution
given by $z^* = z$, $x^*=-y$. A $C^*$-algebra
$C(S_{q,s}^2)$ corresponding to  $A(S_{q,s}^2)$ is  defined if 
$q$ is real, nonzero and $-1<q<1$, and either $s\in [0,1]$ or else 
$s^2 = -q^{2n}$, $n\in {\bf N}$. 
In the case $s^2=-q^{2n}$  the quantum sphere is a finite dimensional
$C^*$-algebra  
isomorphic to the full matrix algebra  
${\rm Mat}_n({\bf C})$, which can  
be interpreted as a deformation of the fuzzy sphere \cite{Mad:fuz}. 
Thus it is  termed a {\em $q$-fuzzy sphere} and 
denoted by $S^2_{q,n}$. $S^2_{q,n}$ 
is no longer a subalgebra of
$SU_q(2)$. On an algebraic level (with $q$ a root of unity), $S^2_{q,n}$
has been shown in \cite{AleRec:non} to describe
$D$-branes of open string in the $SU(2)$
Wess-Zumino-Witten model, and was recently studied in
 \cite{GroMad:fie}. 

The algebraic noncommutative 4-sphere $A(S^4_{p,q,s})$, 
corresponding to the algebra
$A(S_{q,s}^2)$ is a noncommutative algebra generated by 
$1,\xi,\eta,\zeta, U,V$ subject to 
the relations
\begin{eqnarray}
\zeta\xi = q^2\xi\zeta,\quad \eta\zeta = q^2\zeta\eta, \quad \xi U = p U\xi,
\quad V\xi = p\xi V, \nonumber\\
\eta V = pV\eta, \quad U\eta =p\eta U, \quad UV=VU, \quad U\zeta=\zeta
U, \quad V\zeta = \zeta V, \label{rel.s4}
\\\xi\eta =
(\zeta -1)(\zeta +s^2)+UV, \quad \eta\xi = (q^{2}\zeta -1)(q^{2}\zeta
+s^2)+UV, \nonumber
\end{eqnarray}
where $p,q\neq 0$, $s^{2}\neq -1$ are complex parameters. The algebra 
$A(S^4_{p,q,s})$ can be made into a $*$-algebra  with $\zeta^* =
\zeta$, $\xi^* = -\eta$ and $U^*=V$, provided $q^{2}$, $s^2$ are real 
and $p$ is a 
pure phase, i.e., $p = 
\exp(2\pi i \theta)$, where $\theta\in [0,1)$. 
 The corresponding $C^{*}$-algebra is denoted by
$C(S^4_{q,\theta,s})$. Explicitly, $C(S^4_{q,\theta,s})$ 
 is defined by the operator (supremum) norm 
closure over all  admissible $*$-representations of a
polynomial involutive algebra which has a presentation
with generators $\xi,\zeta, U$ and relations
\begin{eqnarray}
\zeta\xi = q^2\xi\zeta, \quad \xi U = p U\xi,
\quad U^*\xi = p\xi U^*, \quad
 UU^*=U^*U, \quad U\zeta=\zeta
U, \nonumber
\\
\xi\xi^* +
(\zeta -1)(\zeta +s^2)+UU^*=0, \quad \xi^*\xi + (q^{2}\zeta -1)(q^{2}\zeta
+s^2)+UU^* = 0, \label{cs4}
\end{eqnarray}
with 
$-1<q<1$, $0< s\leq 1$ and $p = \exp(2\pi i\theta)$. The quantum
4-sphere $S_{q,\theta}^4$ in \cite{Sit:mor} is isomorphic to
$C(S^4_{q^2,\theta,1})$. 
 Unitary
representations $\pi_{c, \pm}$ of $C(S^4_{q,\theta,s})$ in a Hilbert space with basis
$|k,l\rangle$ with $k\in {\bf N}$, $l\in {\bf Z}$ are labelled by 
a complex number $c$, $|c|\leq s$,
and  are given by
$$
\pi_{c, \pm}(U)|k,l\rangle = c|k,l+1\rangle, \qquad 
\pi_{c, \pm}(U^*)|k,l\rangle = c^{*}|k,l-1\rangle,
$$$$ 
\pi_{c, \pm}(\zeta)|k,l\rangle =
\alpha_{\pm}q^{2(k-1)}|k,l\rangle,
$$$$  
\pi_{c, \pm}(\xi)|k,l\rangle =
p^{l}\omega_{\pm,k}|k+1,l\rangle, \qquad 
\pi_{c, \pm}(\xi^*)|k+1,l\rangle =
p^{-l}\omega_{\pm,k}|k,l\rangle,
$$
where 
$$
\alpha_\pm = \frac{1}{2}\left(1-s^2 \pm \sqrt{(s^2+1)^2 -
4|c|^2}\right), \quad \omega_{\pm,k} =
\sqrt{\left(1-\alpha_{\pm}q^{2k}\right)
\left(s^2+\alpha_{\pm}q^{2k}\right)-|c|^{2}},
$$
$k = 0,1,2,\ldots $. Note that if $s=1$ these representations involve representations of
the quantum sphere  $C(S^{2}_{q,1})$ on the Hilbert space spanned by
$|k\rangle$, $k\in {\bf N}$, thus reducing to  
representations of the type discussed in 
\cite{Sit:mor}. 

In the case $s^2 =
-q^{2n}$, the $C^*$-algebra  $C(S^4_{q,\theta,s})$ reduces to 
$C(S^{2}_{q,s})$ and thus is isomorphic to
the full matrix algebra ${\rm Mat}_n({\mathbf C})$. Note, however, that
on a purely algebraic level it does make sense to define
$A(S^4_{q,\theta,s})$ corresponding to the $q$-fuzzy sphere
$S_{q,n}^2$. A study of such a noncommutative space might be of interest 
to string theory.

One can easily construct an example of  a projector associated
to the noncommutative algebra of functions on the quantum 4-sphere 
$A(S^4_{p,q,s})$. 
This defines a noncommutative vector bundle (projective module) 
over $A(S^4_{p,q,s})$. Consider the following $4\times 4$ matrix with 
entries from $A(S^4_{p,q,s})$, 
$$
e = \frac{1}{1+s^{2}}\pmatrix{1-\zeta & 0 &U &\xi \cr
                              0 &1-q^{2}\zeta & -\eta & -pV\cr
			      V &\xi &s^{2}+\zeta & 0\cr
			      -\eta & -p^{-1}U &0 & s^{2}+q^{2}\zeta}.
$$
One easily checks that $e^{2}=e$, hence it is a projector as claimed. 
Furthermore in the case when $A(S^4_{q,\theta,s})$ is a $*$-algebra and
thus including $C(S^4_{q,\theta,s})$, this projector is 
self-adjoint, i.e., $e^{*}=e$ ($*$ in ${\rm 
Mat}_{4}(A(S^4_{q,\theta,s}))$ combines $*$ in $A(S^4_{q,\theta,s})$ 
with matrix transposition). Thus $e$ defines a projective module $E = 
\{{\bf v}e\; | \; {\bf v}= (v_1,v_2,v_3,v_4), \; v_i\in A(S^4_{p,q,s})\}$
or a 
noncommutative vector bundle over $A(S^4_{q,\theta,s})$ \cite{Con:non}. 
Such vector bundles are
classified by the Chern-Connes classes in cyclic homology. 
The components of the 
Chern-Connes class 
of $E$ are defined by
$$
ch_{n}(E) = c_{n}\sum_{i_{1}\ldots i_{2n+1}}(e-\frac{1}{2})_{i_{1}i_{2}}\otimes 
\bar{e}_{i_{2}i_{3}}\otimes \bar{e}_{i_{3}i_{4}}\otimes\cdots\otimes 
\bar{e}_{i_{2n+1}i_{1}},
$$
where $\bar{e}$ is $e$ projected down to the nonunital part of the 
algebra of functions on the quantum sphere $A(S^4_{p,q,s})/{\bf C}1$, 
and $c_{n}$ are normalisation factors. The Chern-Connes class is a 
cocycle in the cyclic homology of $A(S^4_{p,q,s})$. Clearly 
$ch_{0}(E) = 0$ and
\begin{eqnarray*}
ch_{1}(E) &\propto &
\frac{1}{(1+s^{2})^{3}}(q^{2}-1)(\zeta\otimes\left(U\otimes V - 
V\otimes U\right) + U\otimes\left(V\otimes \zeta - \zeta\otimes 
V\right) \\
&&+V\otimes\left(\zeta\otimes U- U\otimes\zeta\right)).
\end{eqnarray*}
Note that up to normalisation and slightly different conventions 
$ch_{1}(E)$ has the same form as the ones computed in 
\cite{DabLan:ins} and \cite{Sit:mor}. 
In particular, within 
the range of $s$ the first factor can be absorbed in the 
normalisation, and thus $ch_{1}(E)$ does not essentially depend on 
$s$. Furthermore, $ch_{1}(E) =0 $ if and only if $q = \pm 1$.

Following the same method as in \cite{Sit:mor} one can define 
a projective module $\tilde{E}$ over $A(S^4_{p,q,s})$ for which 
$ch_{0}(\tilde{E}) = ch_{1}(\tilde{E}) =0$, provided one formally
adjoins a  (self-adjoint) central element $Z = \sqrt{UV}$ to 
$A(S^4_{p,q,s})$. 
In the case of the corresponding 
$C^{*}$-algebra this can be done by considering a suitable infinite 
series.  To construct noncommutative vector bundle
 $\tilde{E}$ one uses the charge 1 magnetic monopole projector for all 
Podle\'s two-spheres constructed in \cite{BrzMaj:lin} 
\cite{BrzMaj:geo},
and defines a projector $\tilde{e}$ in ${\rm Mat}_{4}(A(S^4_{p,q,s}))$ by
\begin{equation}
\tilde{e} = \frac{1}{2(1+s^{2})}\pmatrix{1+s^{2}+2Z & 0 & 
1-s^{2}-2\zeta &2\xi \cr
         0&1+s^{2}+2Z & -2\eta & s^{2}-1 +2q^{2}\zeta\cr
	 1-s^{2}-2\zeta &2\xi &1+s^{2}-2Z & 0\cr
         -2\eta & s^{2}-1+2q^{2}\zeta&0&1+s^{2}-2Z}.
\label{tildee}
\end{equation}
One can easily find directly that $\tilde{e}$ is a projector and that 
$ch_{0}(\tilde{E}) = ch_{1}(\tilde{E}) =0$. There is no need to do it
here, for a general 
justification of this fact is provided below. Clearly 
$\tilde{e}$ is self-adjoint with respect to the $*$-structure on 
$A(S^4_{q,\theta,s})$, whenever defined. Note also that since 
$\tilde{e}$ depends on $Z$, $\zeta$, $\xi$ and $\eta$ (and does not 
depend on $U$ and $V$ separately) it might be viewed as defined on 
the quantum 3-sphere \cite{DabLan:inst}.\medskip

\noindent {\bf 3.} The projector $\tilde{e}$ is a special case of more 
general construction which provides one with a wide range of 
noncommutative algebras and corresponding projectors. Consider an
algebra $A$ generated by $1$, an $n\times n$ matrix of generators ${\bf t} = 
(t_{ij})$ and by an additional generator $Z$. Consider another $n\times 
n$ matrix $\tilde{\bf t}$ of elements of $A$. In case $A$ is a 
${*}$-algebra one requires $Z^{*}=Z$ and $\tilde{\bf t}={\bf 
t}^{*}$. Then a $2n\times 2n$ matrix $e$ with entries from $A$ given 
in the block form by
\begin{equation}
e = \frac{1}{2}\pmatrix{1+Z& {\bf t}\cr \tilde{{\bf t}} & 1-Z},
\label{block}
\end{equation}
is a (self-adjoint) projector provided $Z$ is central in $A$ and 
\begin{equation}
\tilde{\bf t}{\bf t} = {\bf t}\tilde{\bf t} = 1-Z^{2}.
\label{det}
\end{equation}
The first two components of the Chern-Connes character of the 
corresponding noncommutative vector bundle $E$ over $A$ come out as 
$ch_{0}(E) = 0$ and
\begin{eqnarray}
    ch_{1}(E) & \propto & \sum_{ij}(t_{ij}\otimes (
    \tilde{t}_{ji}\otimes Z - Z\otimes \tilde{t}_{ji})+ 
    \tilde{t}_{ji}\otimes (Z\otimes t_{ij}- t_{ij}\otimes Z) \nonumber \\
    &&+ Z\otimes(t_{ij}\otimes \tilde{t}_{ji} - \tilde{t}_{ji}\otimes 
    t_{ij})).\label{ch1}
\end{eqnarray}
This is precisely the method used to obtain $\tilde{e}$ above. As 
proven in \cite{BrzMaj:lin}  \cite{BrzMaj:geo} the following matrix 
with
entries from $A(S^2_{q,s})$,
\[ f = {1\over 1+s^{2}}
\pmatrix{1-z &x\cr -y & s^{2}+q^{2}z},\]
is a projector in ${\rm Mat}_2(A(S^2_{q,s}))$, i.e., 
$f^{2} = f$. The matrix $f$ describes a noncommutative line bundle 
associated to the Dirac monopole principal bundle (Hopf fibration) 
over $A(S^2_{q,s})$. The connection defined by $f$ (the Grassmann 
connection) is a gauge field of the q-deformed Dirac monopole. The 
fact that $f^{2}=f$ implies that 
$$
{\bf t} = \tilde{\bf t} = {2\over 1+s^{2}}
\pmatrix{1-\zeta &\xi\cr -\eta & s^{2}+q^{2}\zeta} -1
$$ 
satisfies ${\bf t}^2 = 1- (2Z/(1+s^2))^2$, where $Z^{2}=UV$ 
in (suitably extended)  $A(S^4_{p,q,s})$, i.e., the condition
(\ref{det}) holds. Matrix $\tilde{e}$ in (\ref{tildee}) has 
precisely the block form (\ref{block}) (with $Z$ replaced by 
$2Z/(1+s^2)$). Since ${\bf t} = \tilde{\bf 
t}$ the first component of the corresponding Chern-Connes character 
vanishes by (\ref{ch1}). In fact again using that ${\bf t} = \tilde{\bf 
t}$ one easily finds that $ch_k(\tilde{E}) =0$, k$=0,1,2,\ldots$

A rich source of algebras with projectors of the block form
(\ref{block}) and non-trivial Chern-Connes characters is 
provided by the standard, FRT-construction of matrix quantum groups. 
Recall that the FRT-construction \cite{FadRes:Lie} associates an 
algebra $A(R)$ to any invertible $n^{2}\times n^{2}$ solution $R$ of the quantum 
Yang-Baxter equation $R_{12}R_{13}R_{23}= R_{23}R_{13}R_{12}$. The 
algebra $A(R)$ is generated by $1$ and an $n\times n$ matrix $\bf t$ 
subject to the following
RTT-relations
\begin{equation}
    R{\bf t}_{1}{\bf t}_{2}= {\bf t}_{2} {\bf t}_{1}R,
\label{frt}
\end{equation}
where ${\bf t}_{1}={\bf t}\otimes 1$, ${\bf t}_{2}=1\otimes {\bf t}$. 
In each of such 
algebras there is a central element known as the {\em quantum 
determinant}, $\det_{q}({\bf t})$ ($\det_{q}({\bf t})$ is grouplike if
$A(R)$ is equipped with the matrix bialgebra structure). 
Quantum determinant can be computed 
explicitly in the framework of braided groups (cf.\ \cite{Maj:book}). 
Furthermore, by using quantum minors in ${\bf t}$ one can construct a 
matrix $\tilde{\bf t}$ with the property ${\bf t}\tilde{\bf t} = 
\tilde{\bf t}{\bf t} = \det_{q}({\bf t})$. Now one can define an 
algebra $A$  by adjoining a central element $Z$ to $A(R)$ which 
is required to satisfy $\det_{q}({\bf 
t}) = 1- Z^{2}$. The resulting algebra $A$ has an associated 
noncommutative vector bundle over itself with a projector given by 
(\ref{block}).  

As an example  consider  an algebra associated to the standard 
solution of the quantum Yang-Baxter equation corresponding to 
$GL_{q}(n)$,
$$
R = q^{-1/n}\left(q\sum_{i}E_{ii}\otimes E_{ii} +\sum_{i\neq j}E_{ii}\otimes 
E_{jj} + (q-q^{-1})\sum_{i>j}E_{ij}\otimes E_{ji}\right),
$$
where $E_{ij}$ are the usual matrix units. The elements of  matrix 
$\tilde{\bf t}$ come out as
$$
\tilde{t}_{ij} = (-q)^{i-j}\sum_{\sigma\in 
S_{n-1}}(-q)^{\ell(\sigma)}t_{j_{1}i_{\sigma(1)}}\cdots 
t_{j_{n-1}i_{\sigma(n-1)}},
$$
where $S_{n-1}$ is the permutation group and $\{i_{1}, \ldots , 
i_{n-1}\} = \{1,\dots, i-1,i+1,\ldots n\}$, and 
$\{j_{1}, \ldots , 
j_{n-1}\} = \{1,\dots, j-1,j+i,\ldots n\}$. In the case $q$-real  one can 
define a consistent $*$-structure by imposing $t_{ij}^{*} = 
\tilde{t}_{ji}$. The resulting algebra $A$ is a ${*}$-algebra given by
 the 
RTT-relations (\ref{frt}) and $t_{ij}^{*}t_{jk} = t_{ij}t_{jk}^{*} = 
(1-Z^{2})\delta_{ik}$, and has associated self-adjoint projector $e$ in
${\rm Mat}_{2n}(A)$ of the form in equation~(\ref{block}). 
In the case $n=2$ this is exactly the quantum 
4-sphere $S_{q}^{4}$ and $e$ the corresponding projector 
introduced in \cite{DabLan:ins}. \medskip

{\bf 4.} In this note we have extended some results of recent papers 
\cite{DabLan:ins} and \cite{Sit:mor}. We introduced a wide class of
examples of noncommutative spaces with implicit quantum group symmetry.
It is  hoped that these examples will help in deciding how axioms for
a noncommutative manifold in \cite{Con:gra} might be modified in order
to include examples based on quantum groups. On the other hand it would
be interesting, and we belive indeed desired, to study whether (some of)
the introduced vector bundles can be viewed as bundles associated to
quantum (coalgebra) principal bundles \cite{BrzMaj:coa} in the way analogous to the
$q$-deformed Dirac monopole  in \cite{BrzMaj:geo} or a
bundle over the quantum 4-sphere in \cite{BonCic:ins}, and 
whether the constructed projectors correspond to strong connections 
\cite{Haj:str}\cite{DabGro:str} on such principal bundles.\bigskip

{\sc Acknowledgements.}
This research is supported by the
British Council grant WAR/992/147. We would like to thank Ludwik D\c 
abrowski and Giovanni Landi for useful comments. T.\ Brzezi\'nski thanks
EPSRC for an Advanced Research Fellowship.



\begin{thebibliography}{99}
\bibitem{AleRec:non} A.Yu.\ Alekseev, A.\ Recknagel and V.\ Schomerus.
Non-commutative world-volume geometries: Branes on $SU(2)$ and fuzzy
spheres. {\em JHEP} 9909, 023 (1999).
\vspace{-3mm}\bibitem{BonCic:ins} F.\ Bonechi, N.\ Ciccoli and M.\ Tarlini. 
Noncommutative instantons and the 4-sphere from quantum groups. {\em 
Preprint} math.QA/0012236. 
\vspace{-3mm}\bibitem{Brz:hom} T.\ Brzezi\'nski. Quantum homogeneous spaces as 
quantum quotient spaces. {\em J.\ Math.\ Phys.} {\bf 37} (1996), 
2388--2399.
\vspace{-3mm}\bibitem{BrzMaj:coa} T.\ Brzezi\'nski and S.\ Majid. Coalgebra bundles. 
{\em Commun.\ Math.\ 
Phys.} {\bf 191} (1998), 467--492.
\vspace{-3mm}\bibitem{BrzMaj:lin} T.\ Brzezi\'nski and S.\ Majid. Line bundles on 
quantum spheres. [in:] {\em Particles, Fields and Gravitation}, J.\ 
Rembieli\'nski (ed.), AIP Woodbury, New York, pp.\ 3--8, 1998.
\vspace{-3mm}\bibitem{BrzMaj:geo} T.\ Brzezi\'nski and S.\ Majid. Quantum geometry 
of algebra factorisations and coalgebra bundles. {\em Commun.\ Math.\ 
Phys.} {\bf 213} (2000), 491--521.
\vspace{-3mm}\bibitem{Con:non} A.\ Connes. {\em Noncommutative Geometry.} Academic 
Press, 1994.
\vspace{-3mm}\bibitem{Con:gra} A.\ Connes. Gravity coupled with matter and
foundations of noncommutative geometry. {\em Commun.\ Math.\ Phys.} {\bf
182} (1996), 155-176.
\vspace{-3mm}\bibitem{ConLan:man} A.\ Connes and G.\ Landi. Noncommutative 
manifolds, the instanton algebra and isospectral deformations. {\em 
Preprint} math.QA/0011194
\vspace{-3mm}\bibitem{DabGro:str} L.\ D\c abrowski, H.\ Grosse and 
P.M.\ Hajac. Strong connections and Chern-Connes pairing in the 
Hopf-Galois theory. {\em Preprint} math.QA/9912239. To appear in {\em Commun.\ Math.\ 
Phys.}
\vspace{-3mm}\bibitem{DabLan:inst} L.\ D\c abrowski and G.\ Landi.
Instanton algebras and quantum 4-spheres. {\em Preprint} math.QA/0101177.
\vspace{-3mm}\bibitem{DabLan:ins} L.\ D\c abrowski, G.\ Landi and T.\ Masuda. 
Instantons on the quantum 4-spheres $S^{4}_{q}$. {\em Preprint} 
math.OA/0012103.
\vspace{-3mm}\bibitem{FadRes:Lie} L.D.\ Faddeev, N.Yu.\ Reshetikhin and L.A.\ 
Takhtajan. Quantisation of Lie groups and Lie algebras {\em Leningrad 
Math.\ J.} {\bf 1} (1990), 193--225.
\vspace{-3mm}\bibitem{GroMad:fie} H.\ Grosse, J.\ Madore and H.\ Steinacker. Field
theory on the $q$-deformed fuzzy sphere I. {\em Preprint}
hep-th/0005273. To appear in {\em J.\ Geom.\ Phys.}
\vspace{-3mm}\bibitem{Haj:str} P.M.\ Hajac. Strong connections on 
quantum principal bundles. {\em Commun.\ Math.\ 
Phys.} {\bf 182} (1996), 579--617.
\vspace{-3mm}\bibitem{Mad:fuz}
J.\ Madore. The fuzzy sphere. {\em Class.\ Quant.\ Grav.} {\bf 9}
(1992), 69--87.
\vspace{-3mm}\bibitem{Maj:book}
S.~Majid. {\em Foundations of Quantum Group Theory}, Cambridge University 
Press 1995.
\vspace{-3mm}\bibitem{Pod:sph} P.\ Podle\'s. Quantum spheres. {\em Lett.\ Math.\
Phys.} {\bf 14} (1987), 193--202.
\vspace{-3mm}\bibitem{Sit:mor} A.\ Sitarz. More noncommutative 4-spheres. {\em 
Preprint} math-ph/0101001.
\end{thebibliography}
\end{document}